
\documentstyle{amsppt}

\TagsOnRight

\nopagenumbers

\input amstex
\magnification=\magstep1
\NoBlackBoxes


                
\def\bet{\beta}         
                \def\Gam{\Gamma}

\def\tet{\theta}

\def\lam{\lambda}

                \def\Ome{\Omega}

\def\F{\Cal F}

\def\fo{\F{|_Z}\otimes o(Z)}

\def\E{\Cal E}
\def\End{\operatorname{End\,}}
\def\Hom{\operatorname{Hom\,}}
\def\ind{\operatorname{ind}}
\def\im{\operatorname{im}}

\def\calM{{\Cal M}}
\def\calN{{\Cal N}}
\def\calP{{\Cal P}}
\def\calQ{{\Cal Q}}

\def\calL{{\Cal L}}

     \def\grg{{\frak g}}
  \def\RR{\Bbb R}  
\def\ZZ{\Bbb Z}  \def\CC{\Bbb C}

\def\w{\tet}

\def\n{\nabla}

\def\Ss{\Cal S}

\def\int{\operatorname {int}}



\topmatter

\title Novikov inequalities with symmetry
\endtitle
\author  Maxim Braverman and Michael Farber \endauthor
\address
School of Mathematical Sciences,
Tel-Aviv University,
Ramat-Aviv 69978, Israel
\endaddress
\email farber\@math.tau.ac.il, maxim\@math.tau.ac.il 
\endemail
\thanks{The research was  supported by grant No. 449/94-1 from the
Israel Academy of Sciences and Humanities}
\endthanks
\abstract
We suggest Novikov type inequalities in the situation of a compact Lie
groups action assuming that the given closed 1-form is invariant and
basic. Our inequalities use equivariant cohomology and an appropriate
equivariant generalization of the Novikov numbers.  We test and apply
our inequalities in the case of a finite group. As an application we
obtain Novikov type inequalities for a manifold with boundary.
\endabstract 
\endtopmatter

In 1981 S.Novikov found a generalization of the classical Morse
inequalities to closed 1-forms.  In this paper we suggest an
equivariant version of the Novikov inequalities.  We consider compact
$G$ manifold $M$, where $G$ is a compact Lie group, and an invariant
closed 1-form $\theta$ on $M$. We assume that the form $\theta$ is
non-degenerate in the sense of Bott, and our problem is to find
estimates on the topology of the set $C$ of critical points of
$\theta$ using global topological invariants of $M$.

We construct {\it the equivariant Morse counting series}, which
combines information about the equivariant cohomology of all the
connected components of $C$. Assuming that the form $\theta$ is {\it
basic} (cf. below) we define an equivariant generalization of the
Novikov numbers and, using these numbers, we construct {\it the
Novikov counting series}.  Our main theorem (Theorem 7) states that
the equivariant Morse series {\it is greater} (in an appropriate
sense) than the Novikov counting series. This statement contains an
infinite number of inequalities involving the dimensions of the
equivariant cohomology of connected components of $C$ and the
equivariant Novikov numbers.

We use in this paper {\it equivariant cohomology twisted by
equivariant flat vector bundles}, which is crucial for our
approach. On one hand, any closed invariant basic 1-form determines a
one-parameter family of equivariant flat vector bundles, which we use
to define the equivariant generalizations of the Novikov numbers.  On
the other hand we observe, that using this cohomology allows to
strengthen the inequalities considerably (cf. \cite{BF1,\S 1.7}.

Simple examples show, that applying the well-known equivariant Morse
inequalities of Atiyah and Bott \cite{AB} to the case when the group
$G$ is finite, one obtains the estimates, which are sometimes worse
than the standard Morse inequalities (ignoring the group action!). The
situation may be improved, however, by using the twisted equivariant
cohomology. If $G$ is a finite group, then any representation of $G$
gives rise to an equivariant flat vector bundle and then (applying our
general construction) to a family of inequalities. Examples show that
only all these inequalities (corresponding to all the irreducible
representations) put together give good enough estimate of the
topology of $C$.

As a simple application we obtain
Novikov type inequalities for manifolds with boundary.


\subheading{1. Basic 1-forms} Consider a closed 1-form $\w$ on a
closed manifold $M$. Our problem is to find estimates on the topology
of the set $C$ of critical points of $\w$ by using global topological
invariants of $M$.  Suppose that $G$ is a compact Lie group acting on
$M$. We will assume that $\w$ is $G$ {\it invariant}, i.e. $g^*\w=\w$
for any $g\in G$. Moreover, we will assume that $\w$ is {\it basic};
recall, that this means that $\w$ vanishes on vectors tangent to the
orbits of $G$.

Note that  any {\it exact} invariant form $\theta=df$ is
basic. Also, (cf., for example, \cite{6,~Lemma~3.4}), 
{\it if $M$ is connected and if the set of fixed points of the action of 
$G$ on $M$ is not empty, then any closed $G$-invariant 1-form on $M$ 
is basic}.

\subheading{2. Equivariant flat vector bundles} Let $\F\to M$ be a
flat vector bundle over $M$ endowed with a smooth action of $G$ such
that the projection $\F\to M$ is $G$-equivariant.  We will assume
that the action $g: \F_x\to \F_{g\cdot x}$ is linear for any $g\in G$
and $x\in M$.

The group $G$ acts naturally on the space of differentials forms
$\Ome^*(M,\F)$ on $M$ with values in $\F$.  For any element $X\in
\grg$, (where $\grg$ denotes the Lie algebra of $G$) we will denote by
$\calL^\F(X):\Ome^\ast(M,\F)\to \Ome^\ast(M,\F)$ the corresponding Lie
algebra action.

A flat bundle $\F\to M$ as above, is called {\it $G$-equivariant flat
vector bundle} if the following two conditions are satisfied: 
$$
	g\circ \n\ =\ \n\circ g, \qquad \n_{X_M}\ =\ \calL^\F(X)
			\tag1 
$$
for any $g\in G$, and for any $X\in \grg$.  Here $\n$ denotes the
covariant derivative of $\F$ and for $X\in \grg$, $X_M$ denotes the
vector field on $M$ defined by the action of $\grg$.  The second
condition determines the covariant derivative in the directions
tangent to the orbits.

If the group $G$ is connected, then first condition
in (1) follows from the second. Conversely, if $G$ is finite,
the second condition in (1) carries no information.  
However, these conditions are independent in general. 

As an {\it example}, consider a closed $G$-invariant 1-form $\w$ on
$M$ with complex values.  Consider the flat vector bundle determined
by the form $\w$. Namely, let $\E_\w =M\times \CC$ with the $G$ action
coming from the factor $M$ and with the flat connection
$\n=d+\w\wedge\cdot$.  This flat bundle always satisfy the first
condition of (1) and it satisfies the second condition iff the form
$\theta$ is basic.

\subheading{3. The pushforward} 
Suppose (only in this section) that {\it
the action of $G$ on $M$ is free}. Then the quotient $B=M/G$ is a 
smooth manifold and the map $q:M\to B$ is a locally trivial
fibration. We will discuss a construction which produces a flat vector 
bundle over $B$ starting from an equivariant vector bundle over $M$.

Namely, let $\Ss(\F)$ denote the locally constant sheaf of flat
sections of $\F$.  Then the direct image $q_*\Ss(\F)$ of $\F$ is a
locally constant sheaf over $B$. Let $q_*\F$ denote the flat bundle
corresponding to this sheaf.  The group $G$ acts naturally on $q_*\F$
and this action is compatible with the flat structure. Moreover, the
action of the connected component $G_0$ of the unit element $e\in G$
is trivial. Thus, the bundle $q_*\F$ splits into a direct sum of its
flat subbundles corresponding to different irreducible representations
of the finite group $G/G_0$. The most important for us will be the
subbundle corresponding to the trivial representation; we will denote
it by $q^G_*\F$. Note that if $G$ is connected, then $q_*^G\F=q_*\F$.
{\it We will say that the flat bundle $q_*^G\F$ is the pushforward of
the bundle $\F$.}

\subheading{4. Twisted equivariant cohomology} Next we define
equivariant cohomology $H^\ast_G(M,\F)$ of $M$ with coefficients in an
equivariant flat vector bundle $\F$.  The idea of the construction is
as follows. Let $EG\to BG$ be the universal principal bundle. Given an
equivariant flat vector bundle $\F$ over $M$, it induced equivariant
flat bundle $p^\ast\F$ over $EG\times M$, where $p:EG\times M\to M$ is
the projection. Now, we want to form the pushforward $q^G_\ast
p^\ast\F$, where $q: EG\times M\to EG\times_G M =M_G$ is the
projection with respect to the diagonal action. The result is a flat
vector bundle over the Borel's quotient $M_G$. Now we define the
equivariant cohomology $H^\ast_G(M,\F)$ as the cohomology of $M_G$
twisted by the flat vector bundle $q^G_\ast p^\ast\F$.

We cannot literally apply the construction of the previous paragraph
since our category is the category of smooth finite dimensional
manifolds and the universal principal bundle $EG\to BG$ is usually
infinite dimensional. The problem may be overcome by using {\it finite
dimensional approximations of $EG$}. We refer to \cite{6} (see also
\cite{2}) for details.

\subheading{5. Equivariant generalization of the Novikov numbers}
Given an equivariant flat bundle $\F$ over $M$ and a closed basic
1-form $\theta$ on $M$ with real values, consider the one-parameter
family $\F\otimes \E_{t\theta}$ of equivariant flat bundles, where
$t\in \RR$, ({\it the Novikov deformation}).  Here $\E_\w$ denotes the
equivariant flat bundle corresponding $\w$, cf. \S 2.  For a fixed $i$
consider the twisted equivariant cohomology
$$
	H^i_G(M,\F\otimes\E_{t\theta}), 
			\qquad\text{where}\quad t\in\RR,\tag2
$$
as a function of $t\in \RR$.  Then \cite{6, Lemma 1.3} there exists
a {\it finite} subset $S\subset\RR$ such that the dimension of the
cohomology $H^i_G(M,\F\otimes\E_{t\theta})$ is constant for $t\notin
S$ and the dimension of the cohomology
$H^i_G(M,\F\otimes\E_{t\theta})$ jumps up for $t\in S$.  The subset
$S$, is called the {\it set of jump points}; the value of the
dimension of $H^i_G(M,\F\otimes\E_{t\theta})$ for $t\notin S$ is
called {\it the background value of the dimension}.

\proclaim{Definition} The $i$-dimensional {\it equivariant Novikov 
number} $\bet_i^G(\xi,\F)$ is defined as the background value of the
dimension of the cohomology $H^i_G(M,\F\otimes\E_{t\theta})$, $t\in\RR.$
\endproclaim

Here $\xi$ denotes the cohomology class of $\w$. Note that a real
cohomology class $\xi\in H^\ast(M,\RR)$ lies in the image of the
natural map $H_G^\ast(M,\RR)\to H^1(M,\RR)$ if and only if it may be
represented by a basic differential form.  Also the equivariant flat
bundle $\E_\theta$ is determined (up to gauge equivalence) only by
$\xi$. Thus, the equivariant Novikov numbers $\bet_i^G(\xi,\F)$ are
well defined for all classes in the image $\xi\in\im[H_G^1(M,\RR)\to
H^1(M,\RR)]$.

The formal power series
$$
	\calN_{\xi,\F}^G(\lam)=
		\sum_i\lam^i\bet_i^G(\xi,\F)\tag3 
$$ 
will be  called {\it the equivariant Novikov series}.

\subheading{6. The equivariant Morse series} Let $C$ denote the set of
critical points of $\w$ (i.e. the set of points of $M$, where $\w$
vanishes).  We assume that $\w$ is {\it non-degenerate in the sense of
Bott}, i.e.  $C$ is a submanifold of $M$ and the Hessian of $\w$ is a
non-degenerate on the normal bundle $\nu(C)$ to $C$ in $M$.  Here by
the Hessian of $\w$ we understand the Hessian of the unique function
$f$ defined in a neighborhood of $C$ and such that $df=\w$ and
$f{|_C}=0$.

Let $Z$ be a connected component of the critical point set $C$ and let
$\nu(Z)$ denote the normal bundle to $Z$ in $M$.  The bundle $\nu(Z)$
splits into the Whitney sum of two subbundles $\nu(Z) = \nu^+(Z)
\oplus \nu^-(Z)$, such that the Hessian is strictly positive on
$\nu^+(Z)$ and strictly negative on $\nu^-(Z)$. The dimension of the
bundle $\nu^-(Z)$ is called the {\it index} of $Z$ (as a critical
submanifold of $\w$) and is denoted by $\ind(Z)$.  Let $o(Z)$ denote
the {\it orientation bundle of $\nu^-(Z)$, considered as a flat line
bundle}.

If the group $G$ is connected, then $Z$ is a $G$-invariant submanifold
of $M$. In general, we denote by $G_Z=\{g\in G| \ g\cdot Z\i Z\}$ the
stabilizer of the component $Z$ in $G$.  Let $|G:G_Z|$ denote the
index of $G_Z$ as a subgroup of $G$; it is always finite.

The compact Lie group $G_Z$ acts on the manifold $Z$ and the flat
vector bundles $\F{|_Z}$ and $o(Z)$ are $G_Z$-equivariant.  Let
$H_{G_Z}^\ast(Z,\fo)$ denote the {\it equivariant cohomology} of the
flat $G_Z$-equivariant vector bundle $\fo$.  Consider the {\it
equivariant Poincar\'e series} 
$$
	\calP^{G_Z}_{Z,\F}(\lam)\ =\ 
	   \sum_i \lam^i \dim_{\CC}H_{G_Z}^i(Z,\fo)\tag4  
$$
and define using it the following {\it equivariant Morse counting
series}
$$
	\calM_{\w,\F}^G(\lam)\ = \ \sum_Z 
	     \lam^{\ind(Z)} |G:G_Z|^{-1}\calP_{Z,\F}^{G_Z}(\lam)\tag5
$$
where the sum is taken over all connected components $Z$ of $C$.

\proclaim{7. Theorem} 
   Suppose that $G$ is a compact  Lie group, acting on a closed manifold
   $M$ and let $\F$ be an equivariant flat vector bundle over $M$. 
   Then for any closed non-degenerate (in the sense of Bott) basic 1-form
   $\w$ on $M$,  there exists a formal power series $\calQ(\lam)$ with
   non-negative integer coefficients, such that
   $$
	\calM_{\w,\F}^G(\lam)- \calN_{\xi,\F}^G(\lam)=
		(1+\lam)\calQ(\lam),\quad\text{where}\quad \xi=[\w].\tag6
   $$
\endproclaim
The proof (cf. \cite{6}) is based in its main
part on the Novikov type inequalities for forms with
non-isolated zeros, obtained in \cite{4,5}.

Consider the case when $G$ acts freely on $M$. Then the basic form
$\theta$ defines a closed 1-form $\theta'$ on $M/G$. It is
straightforward to see, that in this case the inequalities of Theorem
7 (with $\F$ being the trivial line bundle) reduce to the usual
Novikov inequalities with respect to the form $\theta'$ on the
quotient manifold $M/G$.  In particular, we see that in this situation
the equivariant Novikov numbers $\beta_i^G(\xi, \F)$ vanish for large
$i$.

Note that in the case of a (non free) circle action $G=S^1$, the
equivariant Novikov numbers $\beta_i^G(\xi, \F)$ for large $i$ are
two-periodic and coincide with the sum of even or odd (depending on
the parity of $i$) usual (i.e. non-equivariant) Novikov numbers of the
fixed point set; to see this one uses the localization theorem.

An application of Theorem 7 for symplectic torus actions is given in
\cite{6}.

\subheading{8. The finite group case} The rest of this paper is
devoted to illustrations of Theorem 7 in the case when the group $G$
{\it is finite}.  In this situation, one can explicitly calculate the
equivariant cohomology in terms of the action of $G$ on the usual
cohomology.

Let $\rho:G\to \End V_\rho$ be a representation of $G$ on a finite
dimensional complex vector space $V_\rho$. Consider the bundle
$\F_\rho=M\times V_\rho$ over $M$ with the trivial connection and with
the diagonal $G$ action. It is an equivariant flat bundle over $M$,
which is trivial as the flat bundle but it may be not trivial as an
equivariant flat bundle. We will apply Theorem 7 with the bundle
$\F=\F_\rho$.

The twisted equivariant cohomology
$H^*_G(M,\F_\rho\otimes\E_{t\w})$ can be calculated as
$$
	H^*_G(M,\F_\rho\otimes\E_{t\w})=
		\Hom_G\big(V^*_\rho,H^*(M,\E_{t\w})\big).  \tag7
$$
Here $V^*_\rho$ denotes the representation dual to $\rho$; the
cohomology $H^*(M,\E_{t\w})$ is considered with its induced
$G$-action.

For an irreducible representation $\rho$ the equivariant Novikov
number $\beta_i^G(\xi,\F_\rho) = \beta_i^G(\xi,\rho)$ equals to the
{\it background multiplicity} (i.e. generic with respect to $t$) of
$V^\ast_\rho$ in the decomposition of $H^i(M,\E_{t\w})$. The
equivariant Poincar\'e series (4) can be calculated using a formula
similar to (7) in terms of the action of $G$ on the usual cohomology.
We will denote $\calM_{\w,\F_\rho}^G(\lam) = \calM_{\w}^G(\lam;\rho)$
and $\calN_{\xi,\F_\rho}^G(\lam) = \calN_{\xi}^G(\lam;\rho)$ and view
them as functions of two variables: $\lambda$ (which will be formal)
and $\rho$ (which will run over the set of irreducible representations
of $G$).  Applying Theorem 7, we obtain
$$
	\calM_{\w}^G(\lam;\rho)- \calN_{\xi}^G(\lam;\rho)=
		(1+\lam)\calQ(\lam;\rho)  \tag8
$$
where $\rho$ is an irreducible representations of $G$ and 
$\calQ(\lam;\rho)$ denotes a polynomial in $\lambda$ with non-negative 
integral coefficients for any $\rho$. We may view the last statement as 
establishing one family of inequalities of Novikov type for any irreducible 
representation $\rho$.

The inequalities of M.Atiyah and R.Bott \cite{1} correspond (in the case
$\xi=0$) to the inequalities (8) with $\rho$ the trivial representation.
The usual (non-equivariant) Novikov inequalities \cite{7,8} correspond
to the regular representation in (8). Using these remarks, one may construct 
very simple examples (one and two dimensional!) such that the 
non-equivariant inequalities are better than the inequalities of \cite{1}. 

\subheading{9. Novikov inequalities for manifolds with boundary} Our
strategy will be to reduce the problem on a manifold with boundary to
a problem on the double with its natural $\ZZ_2$ action. From section
8 we know that we should expect two families of inequalities -- since
there are two irreducible representations of $\ZZ_2$.

Let $M$ be a compact manifold with boundary $\Gam=\partial M$ and let
$\theta$ be a closed 1-form on $M$. We will denote by $C$ the set of
all critical points of $\theta$. We will suppose that
$\theta|_{\int(M)}$ and $\theta|_\Gamma$ are both {\it non-degenerate in
the sense of Bott} (cf. \S 6) and that {\it any critical point of
$\theta|_\Gamma$ is a critical point of $\theta$ as well}.
Additionally, we will suppose that for any connected component
$Z\subset \Gamma$ of the critical point set of $\theta|_\Gamma$ holds
either 

\roster
  \item {\it $Z$ is nondegenerate as a critical manifold of $\theta$, or
  \item $Z$ is the boundary of a connected component $Z'\subset C$
	such that $Z=Z'\cap\Gamma$ and the intersection 
	$Z'\cap \Gamma$ is transversal.}
\endroster

If the first possibility holds, then the Hessian $h_\tet(\cdot,\cdot)$
of $\theta$ is a nondegenerate quadratic form on the normal bundles to
$Z$ in $M$ and in $\Gamma$ and so there exists a unique nonvanishing
vector field $X$ on $Z$ normal to $\Gamma$ such that $h_\tet$ splits
as a direct sum: $h_\tet(X,Y)=0$ for any $Y\in T\Gamma$.  The
function $h_\tet(X,X)$ is everywhere positive or negative on $Z$; we
will call the corresponding component $Z\subset \Gamma\cap C$ {\it
positive} or {\it negative} respectively.

Represent $C$ as the union of 4 disjoint submanifolds $C_{in}\cup
C_+\cup C_-\cup C_{bd}$, where $C_{in}$ denotes the union of the
connected components of $C$ which do not intersect $\Gamma$, $C_{bd}$
denotes the union of the components which are manifolds with nonempty
boundary, and $C_\pm$ denotes the union of the positive (negative)
components in $\Gamma$ (as defined in the previous paragraph).

For simplicity we will assume that $M$ and all submanifolds $Z\subset
C$ are {\it orientable }-- without this assumption the notations will
be more complicated.

For any connected component $Z\subset C$ denote by $\ind_+(Z)$ and
$\ind_-(Z)$ the dimensions of the positive and the negative subbundles
$\nu_+(Z)$ and $\nu_-(Z)$ of the normal bundle $\nu(Z)$ in $M$
correspondingly, compare \S 6.  Now we will define two {\it Morse
counting polynomials}
$$
	\calM^\pm_\theta(\lambda) = 
		\sum _{Z\subset C_{in}\cup C_{bd}\cup C_{\pm}}
		\lambda^{\ind_\pm(Z)} \calP_Z(\lambda),   \tag9
$$
where $\calP_Z(\lambda)=\sum\lambda^i\dim H^i(Z, o(Z))$ is the Poincar\'e 
polynomial of $Z$; in this formula $Z$ runs over the connected components
contained in $C_{in}\cup C_{bd}\cup C_{+}$ in the case of $+$ and contained 
in $C_{in}\cup C_{bd}\cup C_{-}$ in the case of $-$.

{\it The Novikov numbers} $\beta_i(\xi)$ will be defined as {\it the
background values of the dimension} of $H^i(M, \E_{t\theta})$, cf. 
\S 5. Here $\xi\in H^1(M,\RR)$ denotes the class of $\theta$.

\proclaim{10. Theorem} 
  Under the assumption described above holds
  $$
	\calN_\xi(\lambda) - \calM_\theta^\pm(\lambda) = 
		(1 + \lambda)\calQ^\pm(\lambda),
			\quad\text{where}\quad 
	\calN_\xi(\lambda) = \sum\lambda^i\beta_i(\xi),  \tag10
  $$
  where $\calQ^\pm(\lambda)$ are polynomials with nonnegative integral 
  coefficients.
\endproclaim

Note that in the case of a closed manifold, these two statements (for
$\pm$ equal to $+$ and $-$ correspondingly) are equivalent (as follows
from Poincar\'e duality), but for the case of manifolds with boundary
they are independent.

The proof of Theorem 10 follows by applying Theorem 7 to the double
$D(M)$, two copies of $M$ glued along the boundary $\Gamma$.  The
double $D(M)$ has the natural $\ZZ_2$ action. Using our assumptions on
$\theta$ we may construct an appropriate invariant 1-form $\tilde
\theta$ on $D(M)$ and then use Theorem 7.  The computation of the
twisted equivariant cohomology of the double is based on the following
simple Lemma:

\proclaim{Lemma} Let $\E\to D(M)$ be a flat $\ZZ_2$ equivariant vector
  bundle over the double of a compact manifold with boundary. Suppose
  that the action of $\ZZ_2$ on $\E|_{\partial M}$ is trivial.  Then
  the twisted equivariant cohomology
  $H^i_{\ZZ_2}(D(M),\E\otimes\F_\rho) $ is isomorphic to
  $H^i(M,\E|_M)$, if $\rho$ is the trivial representation, and to
  $H^i(M,\partial M,\E|_M)$, if $\rho$ is the not trivial irreducible
  representation of $\ZZ_2$.  
\endproclaim


\enddocument